\documentclass{aip-cp}

\usepackage[numbers]{natbib}
\usepackage{rotating}
\usepackage{graphicx}
\usepackage{amsfonts}
\usepackage{amssymb}
\usepackage{amstext}
\usepackage{amsmath}
\usepackage{tikz}
\usepackage{wrapfig}

\newcommand{\V}[1]{\mathbf{#1}}
\newcommand{\M}[1]{\mathbf{#1}}

\begin{document}

\title{Scalable Analysis and Design Using Automatic Differentiation}

\author[aff1]{Julian Andrej}
\author[aff1]{Tzanio Kolev}
\author[aff2,cor1]{Boyan Lazarov}

\affil[aff1]{Center for Applied Scientific Computing, Lawrence Livermore National Laboratory, Livermore, CA}
\affil[aff2]{Computational Engineering Division, Lawrence Livermore National Laboratory, Livermore, CA}
\corresp[cor1]{Corresponding author: lazarov2@llnl.gov}

\maketitle

\begin{abstract}
This article aims to demonstrate and discuss the applications of automatic differentiation (AD) for finding derivatives in PDE-constrained optimization problems and Jacobians in non-linear finite element analysis. The main idea is to localize the application of AD at the integration point level by combining it with the so-called Finite Element Operator Decomposition. The proposed methods are computationally effective, scalable, automatic, and non-intrusive, making them ideal for existing serial and parallel solvers and complex multiphysics applications. The performance is demonstrated on large-scale steady-state non-linear scalar problems. The chosen testbed, the MFEM library, is free and open-source finite element discretization library with proven scalability to thousands of parallel processes and state-of-the-art high-order discretization techniques.
\end{abstract}

\section{INTRODUCTION}
Automatic differentiation (AD) \cite{Griewank2008}, or algorithmic differentiation, provides exact values of the Jacobian for complex functions. Despite its long history and many implementations, it remains underutilized in the scientific community. AD simplifies the evaluation of functions into easy-to-differentiate operations, applying the chain rule. Software libraries automate the process, letting researchers focus on their problems rather than differentiating the functions of interest. AD can be implemented in two ways: code transformation and operator/function overloading. Code transformation is based on compiler tools that transform function code into one that evaluates partial derivatives. It requires specific compilers and tools, which may limit platform availability. Modern object-oriented languages like C++ can deploy operator/function overloading, i.e., overload computational operations to compute gradients alongside evaluations. Compared to code transformation, the approach only requires the language compiler without additional tools. In addition, AD can be implemented in two modes: a forward mode, which computes derivatives during function evaluation, reducing memory requirements, and a reverse mode, which involves first evaluating the function, recording operations, and their derivatives. The derivative information propagates through the recorded evaluation tree as a second step. 
Depending on the mode and the implementation, AD adds overhead compared to the standard evaluation process. The overhead can be significant for long and complex computations, and depending on the mode, it can impact the system's memory utilization or the computational cost. Both will significantly affect the total execution time, especially if AD is applied naively on large production codes \cite{Norgaard2017}. Therefore, for finite element (FEM) analysis, we propose to limit the application of AD only to specific parts of the code, preserving the same parallel scalability and performance available to the original code without AD. Furthermore, the proposed approach can be extended to design and optimization problems without any significant coding effort, automating the optimization completely. The proposed localization provides a fast and efficient solution regardless of the AD implementation and the deployed evaluation modes.

\section{AUTOMATIC DIFFERENTIATION IN FINITE ELEMENT ANALYSIS}
\begin{figure*}[th]
    \centering\newcommand{\RoundRect}[4]{
  \draw[
    rounded corners=5,
    black!60!white,
    fill=black!5!white,
  ] (#1,#2) rectangle ++(#3,#4);
}

\newcommand{\DofSquare}[5]{
  \draw[black] (#1,#2) rectangle ++(#3,#4);
  \draw node[fill,circle,inner sep=0pt,minimum size=2.5pt,#5] at (#1, #2) {};
  \draw node[fill,circle,inner sep=0pt,minimum size=2.5pt,#5] at ({#1+#3/2}, #2) {};
  \draw node[fill,circle,inner sep=0pt,minimum size=2.5pt,#5] at ({#1+#3}, #2) {};
  \draw node[fill,circle,inner sep=0pt,minimum size=2.5pt,#5] at (#1, {#2+#4/2}) {};
  \draw node[fill,circle,inner sep=0pt,minimum size=2.5pt,#5] at ({#1+#3/2}, {#2+#4/2}) {};
  \draw node[fill,circle,inner sep=0pt,minimum size=2.5pt,#5] at ({#1+#3}, {#2+#4/2}) {};
  \draw node[fill,circle,inner sep=0pt,minimum size=2.5pt,#5] at (#1, {#2+#4}) {};
  \draw node[fill,circle,inner sep=0pt,minimum size=2.5pt,#5] at ({#1+#3/2}, {#2+#4}) {};
  \draw node[fill,circle,inner sep=0pt,minimum size=2.5pt,#5] at ({#1+#3}, {#2+#4}) {};
}

\newcommand{\QuadSquare}[5]{
  \draw[black] (#1,#2) rectangle ++(#3,#4);
  \draw node[fill,circle,inner sep=0pt,minimum size=2.5pt, #5] at (#1+#3/3, #2+#4/3) {};
  \draw node[fill,circle,inner sep=0pt,minimum size=2.5pt, #5] at (#1+2*#3/3, #2+#4/3) {};
  \draw node[fill,circle,inner sep=0pt,minimum size=2.5pt, #5] at (#1+#3/3, #2+2*#4/3) {};
  \draw node[fill,circle,inner sep=0pt,minimum size=2.5pt, #5] at (#1+2*#3/3, #2+2*#4/3) {};
}

\newcommand{\DofGrid}[5]{%
  \DofSquare{#1}{#2}{#3}{#4}{#5}
  \DofSquare{#1+#3}{#2}{#3}{#4}{#5}
  \DofSquare{#1}{#2+#4}{#3}{#4}{#5}
  \DofSquare{#1+#3}{#2+#4}{#3}{#4}{#5}
}

\def\spacingfactor{0.85}

\newcommand{\DofElem}[5]{%
  \DofSquare{#1}{#2}{#3*\spacingfactor}{#4*\spacingfactor}{#5}
  \DofSquare{#1+#3/\spacingfactor}{#2}{#3*\spacingfactor}{#4*\spacingfactor}{#5}
  \DofSquare{#1}{#2+#4/\spacingfactor}{#3*\spacingfactor}{#4*\spacingfactor}{#5}
  \DofSquare{#1+#3/\spacingfactor}{#2+#4/\spacingfactor}{#3*\spacingfactor}{#4*\spacingfactor}{#5}
}

\newcommand{\QuadElem}[5]{%
  \QuadSquare{#1}{#2}{#3*\spacingfactor}{#4*\spacingfactor}{#5}
  \QuadSquare{#1+#3/\spacingfactor}{#2}{#3*\spacingfactor}{#4*\spacingfactor}{#5}
  \QuadSquare{#1}{#2+#4/\spacingfactor}{#3*\spacingfactor}{#4*\spacingfactor}{#5}
  \QuadSquare{#1+#3/\spacingfactor}{#2+#4/\spacingfactor}{#3*\spacingfactor}{#4*\spacingfactor}{#5}
}

\newcommand{\RoundRectGrid}[2]{
  \RoundRect{0}{0}{#1}{#2}
  \RoundRect{#1 + 0.075}{0}{#1}{#2}
  \RoundRect{0}{#2 + 0.075}{#1}{#2}
  \RoundRect{#1 + 0.075}{#2 + 0.075}{#1}{#2}
}

\begin{tikzpicture}[scale=1.75]

  \RoundRect{0}{0}{0.975}{0.975}
  \RoundRect{1.05}{0}{0.725}{0.975}
  \RoundRect{0}{1.05}{0.975}{0.725}
  \RoundRect{1.05}{1.05}{0.725}{0.725}

  \DofGrid{0.1}{0.1}{0.4}{0.4}{red!75!black}
  \DofGrid{0.9}{0.1}{0.4}{0.4}{red!75!black}
  \DofGrid{0.1}{0.9}{0.4}{0.4}{red!75!black}
  \DofGrid{0.9}{0.9}{0.4}{0.4}{red!75!black}

  \draw[->, line width=0.5pt] (1.85, 0.9) -- node[midway,above] {$P$} ++(.3,0.0);
  \draw[<-, line width=0.5pt] (1.85, 0.8) -- node[midway,below] {$P^T$} ++(.3,0.0);
  \node[align=center] at (0.9, -0.2) {\small T-vector};
  \node[align=center] at (0.9, 2) {\small Global true dofs};

  \begin{scope}[shift={(2.2,0)}]
    \RoundRectGrid{0.85}{0.85}

    \DofGrid{0.1}{0.1}{0.325}{0.325}{black}
    \DofGrid{1.025}{0.1}{0.325}{0.325}{black}
    \DofGrid{0.1}{1.025}{0.325}{0.325}{black}
    \DofGrid{1.025}{1.025}{0.325}{0.325}{black}

    \draw[->, line width=0.5pt] (1.85, 0.9) -- node[midway,above] {$G$} ++(.3,0.0);
    \draw[<-, line width=0.5pt] (1.85, 0.8) -- node[midway,below] {$G^T$} ++(.3,0.0);
    \node[align=center] at (0.9, -0.2) {\small L-vector};
    \node[align=center] at (0.9, 2) {\small Local subdomain dofs};
  \end{scope}

  \begin{scope}[shift={(4.4,0)}]
    \RoundRectGrid{0.85}{0.85}

    \DofElem{0.1}{0.1}{0.325}{0.325}{blue!65!black}
    \DofElem{1.025}{0.1}{0.325}{0.325}{blue!65!black}
    \DofElem{0.1}{1.025}{0.325}{0.325}{blue!65!black}
    \DofElem{1.025}{1.025}{0.325}{0.325}{blue!65!black}

    \draw[->, line width=0.5pt] (1.85, 0.9) -- node[midway,above] {$B$} ++(.3,0.0);
    \draw[<-, line width=0.5pt] (1.85, 0.8) -- node[midway,below] {$B^T$} ++(.3,0.0);
    \node[align=center] at (0.9, -0.2) {\small E-vector};
    \node[align=center] at (0.9, 2) {\small Element dofs};
  \end{scope}

  \begin{scope}[shift={(6.6,0)}]
    \RoundRectGrid{0.85}{0.85}

    \QuadElem{0.1}{0.1}{0.325}{0.325}{green!50!black}
    \QuadElem{1.025}{0.1}{0.325}{0.325}{green!50!black}
    \QuadElem{0.1}{1.025}{0.325}{0.325}{green!50!black}
    \QuadElem{1.025}{1.025}{0.325}{0.325}{green!50!black}

    \draw [->] (1.85, 1.1) to[out=0,in=0,looseness=1.5] node[midway,right] {$D$} (1.85, 0.7) ;
    \node[align=center] at (0.9, -0.2) {\small Q-vector};
    \node[align=center] at (0.9, 2) {\small Quadrature point values};
  \end{scope}


\end{tikzpicture}
    \caption{{ Finite element operators, $A_p$, have a natural decomposition, $A_p = P^T G^T B^T D B G P$, which exposes multi-level parallelism and allows for AD-friendly, matrix-free, memory-efficient implementations that assemble and store only the innermost, pointwise operator component (partial assembly, cf.~\cite{Anderson2021}).}}
    \label{fig_feod}
\end{figure*}
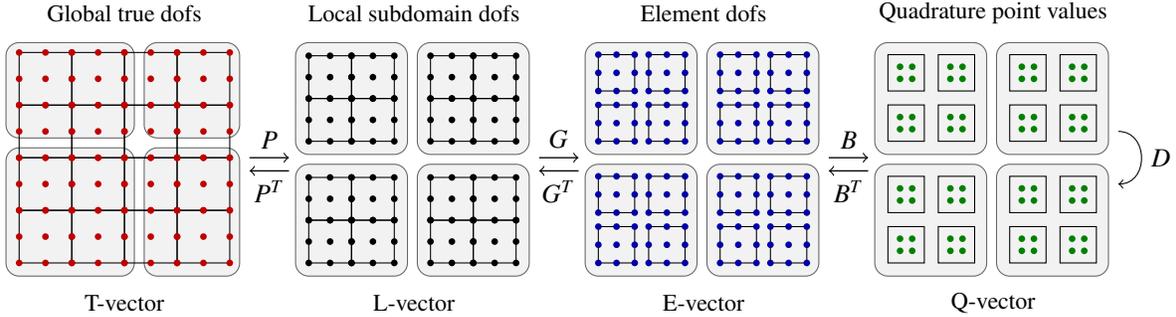

The proposed application of AD in FEM analysis relies heavily on the so-called finite element operator decomposition (FEOD) \cite{Anderson2021} and is demonstrated in Figure \ref{fig_feod}. The subdomain restriction operator $P$ transfers FEM degrees-of-freedom (DOF) from the global to the local subdomain level. The element restriction operator $G$ transfers DOFs from a subdomain level to an element level, and the operator $B$  maps the solution field on the element level to its gradients or values on the quadrature point level. The operator, $D$, is entirely local and is evaluated pointwise at every quadrature point. The decomposition is implemented and available in the MFEM library \cite{Andrej2024}, a free, open-source C++ finite element discretization library. The library is GPU-accelerated with state-of-the-art performance on small user laptops, desktop computer systems, and large high-performance computing (HPC) systems. FEOD encapsulates a generic description of an assembly procedure in a finite element library and allows MFEM to handle derivatives at the innermost level at the quadrature points (D). Operators that transfer data from the global level to subdomain, element, and quadrature levels (P, G, and B) are linear and topological. They do not depend on the solution, physical coordinates, or design parameters and, as a result, are excluded from the differentiation loop, saving both memory and computational resources. The decomposition confines the code modifications to the integration point level, allowing complete automation of the discretization process for complex non-linear problems. The quadrature point-level derivatives can be generated by leveraging Enzyme \cite{Moses2021}, CoDiPack \cite{Sagebaum2019}, or a native MFEM’s internal dual number type implementation. 

For non-linear problems, the finite element operator depends on the solution field $\V{u}$, and its action on $\V{u}$ can be written as
\begin{equation}
\label{eq1}
	A_p\left(\V{u}\right)=P^T G^T B^T D\left(B G P \V{u}\right)\,.
\end{equation}
Differentiating Equation \ref{eq1} results in the following expression for the Jacobian operator
\begin{equation}
\label{eq2}
J_p\left(\V{u}\right)=P^T G^T B^T J_D\left(\V{u}_q\right)B G P \,,
\end{equation}
where $\V{u}_q=B G P\V{u}$ and $J_D\left(\V{u}_q\right)={\text d} D\left(\V{u}_q\right)/ {\text d} \V{u}_q$. The AD application in Equation \ref{eq2} is confined only at the integration point level, i.e., at the constitutive relations, and does not impact the rest of the operators. The latter simplifies the implementation of AD as it requires code modification only at the constitutive relation level. The parallel scalability of the code is not impacted, and any complications arising in black box implementations in parallel are removed. Locally, the constitutive relations can be differentiated using reverse or forward mode. Reverse AD mode introduces additional overhead for memory management and for relations with vector length $\V{u}_q$ equal to the output of the vector function $D\left( \V{u}_q\right)$ forward AD mode will be preferable in terms of computational time.

To demonstrate the advantages of the proposed approach, we have implemented continuous Galerkin finite element discretization of a p-Laplacian problem \cite{Toulopoulos2017} in MFEM. The average computational time and floating point operations (FLOPs) per element for computing the tangent matrix on a cube meshed using 200K elements are reported in Table \ref{tab:a}. The numerical experiments are performed on 12 MPI processes and executed on Intel(R) Xeon(R) CPU E5-2680v4 \@ 2.40 GHz. The timing is obtained using the Caliper library \cite{Boehme2016}, a performance analysis toolbox developed at LLNL, and the results are averaged over 100 runs. The FLOPs are estimated using the Performance Application Programming Interface (PAPI) library \cite{Jagode2019}. The presented results are limited to the assembly of the tangent matrices, and more elaborated analysis and presentation discussing implications for tangent matrix-vector products and residual computations relevant for full automatization of matrix-free non-linear solvers are left for following papers.

\begin{table}[h]
\caption{AD performance for constructing tangent element matrices of tetrahedral elements. RES denotes evaluation using the proposed AD approach, ELM-element level evaluation, and HND-hand coded implementation. Reverse and Forward modes are using CoDiPack, and MFEM denotes the native AD implementation based on dual numbers. }
\label{tab:a}
\tabcolsep7pt\begin{tabular}{|l|r  r|r  r| r  r |r|}
        \hline
        &    \multicolumn{2}{c|}{Reverse} & \multicolumn{2}{c|}{Forward} & \multicolumn{2}{c|}{MFEM} & \\
                    & RES &  ELM & RES &  ELM & RES &  ELM & HND \\
        \hline
	\multicolumn{8}{|c|}{First order tetrahedral element  Tet1 - $\M{K}_e\in \mathbb{R}^{4\times4},\, \V{r}^e \in \mathbb{R}^4$} \\
	\hline
        $\M{K}_e$ [s]     & 0.37 &  0.36 & 0.34 &  0.45 & 0.34 &  0.45 & 0.29 \\
        $\M{K}_e$ [KFLOP] & 3    & 2   & 4    & 10   & 4  & 10 & 2    \\
        \hline
	\multicolumn{8}{|c|}{Second order tetrahedral element  Tet2 - $\M{K}_e\in \mathbb{R}^{10\times10},\, \V{r}^e \in \mathbb{R}^{10}$} \\
        \hline
        $\M{K}_e$ [s]     & 0.85 &  1.57 & 0.81 &  3.40 & 0.80 &  3.31 & 0.62 \\
        $\M{K}_e$ [KFLOP] & 43   &  34  & 45   &  243   & 46 &  242 & 29    \\
        \hline
	\multicolumn{8}{|c|}{Third order tetrahedral element  Tet3 - $\M{K}_e\in \mathbb{R}^{20\times20},\, \V{r}^e \in \mathbb{R}^{20}$} \\
        \hline
        $\M{K}_e$ [s]     & 3.05 &  11.90 & 2.86 &  31.48 & 2.88 &  30.96 & 2.53 \\
        $\M{K}_e$ [KFLOP] & 413   &  388& 419  &  3925  & 424 &  3879 & 279    \\
        \hline
\end{tabular}
\end{table}

AD applications are most commonly found in current finite element literature and codes at the element level instead of the currently proposed integration point level. Denoting the element DOFs vector with $\V{u}_e$ and the element residual with $\V{r}_e$, the element tangent (stiffness) matrix can be expressed as
\begin{equation}
\label{eq3}
	\M{K}_e=\frac{\partial \V{r}_e}{\partial \V{u}_e} = B^T J_D\left(\V{u}_q\right)B\, ,
\end{equation}
where $\V{u}_q=B\V{u}_e$, and $\V{r}_e= B^T D\left(B\V{u}_e\right)$. For the lowest-order linear Lagrangian elements and scalar field problems, like nonlinear diffusion, the number of integration points is relatively low, and the overhead of including the operator $B$ in the differentiation loop is small. However, it's crucial to note that the impact of the operator $B$ becomes significant for high-order elements. Close inspection of Table \ref{tab:a} reveals that the number of floating point operations (FLOP) per element scales proportionally to the number of elemental DOFs compared to the FLOPs required for integration point level AD, with computational time following the same trend. The only exception is the case using reverse AD. Building a computational tree during the forward pass allows for simplifications and optimizations at the cost of more considerable processing (overhead) time. Thus, even though the FLOPs per element have decreased by a factor of 10 for the third-order elements, the average computational time is only reduced by a factor of three. In addition, as discussed and demonstrated in \cite{Norgaard2017}, the reverse mode requires a significant amount of memory for storing the computational tree in contrast to the forward AD mode. Regardless of the reduction in computational cost, the average computational time per finite element is larger than that per forward AD applied at the integration point level. Furthermore, the native MFEM implementation based on dual numbers performs as well as implementations in dedicated libraries, removing the necessity of linking and compiling against external libraries. Finally, the forward AD is only 10-15\% slower than the highly optimized hand-coded implementation of the problem and can be reduced further by deploying more aggressive compilation flags. 

\begin{figure}[ht]
\centering{
\includegraphics[width =.3\textwidth]{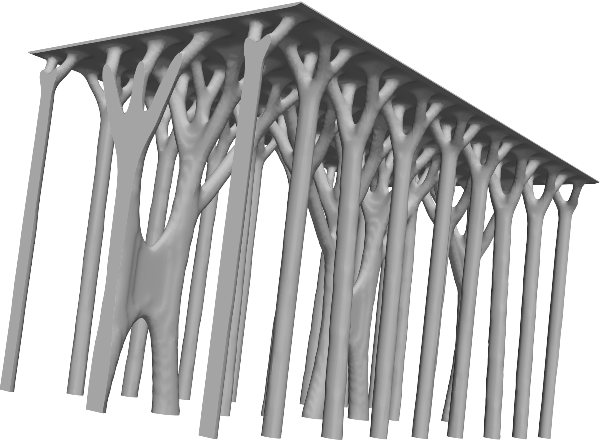}
\includegraphics[width =.3\textwidth]{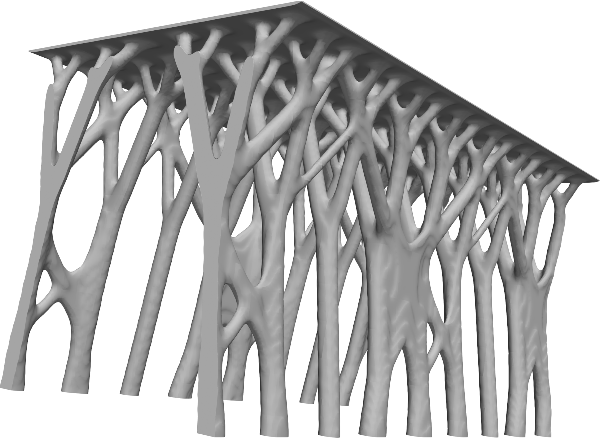}}
\caption{Topology-optimized 3D bridge structure \cite{Duswald2024}.}
\label{fig:fig1}
\end{figure}

The derivatives of any linear or non-linear function (functional in continuous settings), constrained by a discretized PDE written in a residual form as $\mathbf{r}\left(\mathbf{u};\boldsymbol{\rho}\right)=0$, with respect to the parameters ${\boldsymbol{\rho}}$ can be obtained either by adjoint analysis or with the help of AD.  Direct AD   penalizes computational performance in exchange for faster and easier implementation.  The procedure works fine in serial settings and relatively small academic problems. However, in parallel environments, the computational implementation requires the AD tools to account for information exchange between the different processes and, in addition, a considerable amount of memory to accommodate the computational history for realistic simulations. Thus, instead of applying AD as a black-box tool, the suggested technique can be employed only for local computational operations, saving both memory and computational resources. The approach has been demonstrated for Lattice Boltzmann simulations and optimization in \cite{Norgaard2017} and here, we include demonstration  (Figure \ref{fig:fig1}) for topology optimization of large-scale solid mechanics problem \cite{Duswald2024}. The AD is applied on the integration point level for computing adjoint loads and allows for speedy implementation of new objectives and constraints. 

\section{CONCLUSIONS}
Automatic differentiation is a compelling technique applicable to both newly developed applications and existing codes. Careful deployment allows the MFEM library to find derivatives and Jacobians with negligible coding effort. The technique is application-agnostic, and although it is demonstrated here for a steady-state problem, it is applicable to any time-dependent complex multiphysics set of equations.


\section{ACKNOWLEDGMENTS}
This work (LLNL-CONF-866443) was performed under the auspices of the U.S.~Department of Energy by Lawrence Livermore National Laboratory under Contract DE-AC52-07NA27344, the LLNL-LDRD Program under Project tracking No.~22-ERD-009, and Differentiating Large-Scale Finite Element Applications project supported by the U.S.~Department of Energy, Office of Science, Office of Advanced Scientific Computing Research.


\nocite{*}
\bibliographystyle{aipnum-cp}%
\bibliography{bib}%

\end{document}